\documentclass{article}
\usepackage{amssymb}
\usepackage{amsmath}

\input{tcilatex}

\begin{document}

\title{Pseudoprocesses governed by higher-order fractional differential
equations}
\author{Luisa Beghin\thanks{%
Dipartimento di Statistica, Probabilit\`{a} e Statistiche Applicate,
University of Rome ``La Sapienza'', p.le A.Moro 5, 00185 Rome, Italy. E-mail
address: luisa.beghin@uniroma1.it.}}
\maketitle

\begin{abstract}
We study here a heat-type differential equation of order $n$ greater than
two, in the case where the time-derivative is supposed to be fractional. The
corresponding solution can be described as the transition function of a
pseudoprocess $\Psi _{n}$ (coinciding with the one governed by the standard,
non-fractional, equation) with a time argument $\mathcal{T}_{\alpha }$ which
is itself random. The distribution of $\mathcal{T}_{\alpha }$ is presented
together with some features of the solution (such as analytic expressions
for its moments).


\end{abstract}

\section{Introduction}

The study of diffusion equations with a fractional derivative component have
been firstly motivated by the analysis of thermal diffusion in fractal media
in Nigmatullin (1986) and Saichev and Zaslavsky (1997). This topic has been
extensively treated in the probabilistic literature since the end of the
Eighties: see, for examples, Wyss (1986), Schneider and Wyss (1989),
Mainardi (1996), Angulo et al. (2000). Recently fractional equations of
different types have been also studied, such as, for example, the Black and
Scholes equation (see Wyss (2000)) and the fractional diffusion equations
with stochastic initial conditions (see Anh and Leonenko (2000)).

Our aim will concern the extension, to the case of fractional
time-derivative, of a class of equations which is well known in the
literature, namely the \emph{higher-order heat-type equations}. Therefore we
will be interested in the solution of the following problem, for $0<\alpha
\leq 1,\;n\geq 2,$%
\begin{equation}
\left\{
\begin{array}{l}
\frac{\partial ^{\alpha }}{\partial t^{\alpha }}u(x,t)=k_{n}\frac{\partial
^{n}}{\partial x^{n}}u(x,t)\qquad x\in \mathbb{R},\;t>0, \\
u(x,0)=\delta (x)%
\end{array}%
\right.  \label{uno.3}
\end{equation}%
where $\delta (\cdot )$ is the Dirac delta function, $k_{n}=(-1)^{q+1}$ for $%
n=2q,$ $q\in \mathbb{N}$, while $k_{n}=\pm 1$ for $n=2q+1.$ The fractional
derivative appearing in (\ref{uno.3}) is meant, in the Dzherbashyan-Caputo
sense, as%
\begin{equation*}
(D^{\alpha }f)(t)=\frac{d^{\alpha }}{dt^{\alpha }}f(t)=\left\{
\begin{array}{l}
\frac{1}{\Gamma (m-\alpha )}\int_{0}^{t}\frac{f^{(m)}(z)}{(t-z)^{1+\alpha -m}%
}dz\text{,\qquad for }m-1<\alpha <m \\
\frac{d^{m}}{dt^{m}}f(t)\text{,\qquad for }\alpha =m%
\end{array}%
\right. ,
\end{equation*}%
where $m-1=\left\lfloor \alpha \right\rfloor $ and $f\in C^{m}$ (see Samko
et al. (1993) for a general reference on fractional calculus).

In the non-fractional case (which can be obtained from ours as a particular
case, for $\alpha =1$) the \emph{pseudoprocesses} $\Psi _{n}=\Psi _{n}(t),$ $%
t>0$ driven by $n$-th order equations, i.e.
\begin{equation}
\frac{\partial }{\partial t}p(x,t)=k_{n}\frac{\partial ^{n}}{\partial x^{n}}%
p(x,t),\qquad x\in \mathbb{R},\;t>0,  \label{n.1}
\end{equation}%
for $n>2,$ have been introduced in the Sixties and studied so far by many
authors starting from Krylov (1960), Daletsky (1969). Moreover the
distributions of many functionals of $\Psi _{n}$ have been obtained: in
Hochberg and Orsingher (1994) the distribution of sojourn time on the
positive half-line is presented, for $n$ odd, while for an arbitrary $n$ the
same topic is analyzed in Lachal (2003). For $n=3,4,$ the case where the
pseudoprocess is constrained to be zero at the end of the time interval is
considered in Nikitin and Orsingher (2000) and the corresponding
distribution of the sojourn time is evaluated. In Beghin et al. (2000) the
distribution of the maximum is obtained under the same circumstances. In the
unconditional case the maximal distribution is presented in Orsingher
(1991), for $n$ odd, while the joint distribution of the maximum and the
process for diffusion of order $n=3,4$ is presented in Beghin et al. (2001).
Lachal (2003) has extended these results to any order $n>2.$

Some other functionals, such as the first passage time, are treated in
Nishioka (1997) and Lachal (2007). Finally in Beghin and Orsingher (2005) it
is proved that the local time in zero possesses a proper probability
distribution which coincides with the (folded) solution of a fractional
diffusion equation of order $2(n-1)/n,$ $n\geq 2$.

In the fractional case under investigation (i.e. for $0<\alpha <1$) we prove
that the process related to (\ref{uno.3}) is a pseudoprocess $\Psi _{n}$
evaluated at a random time $\mathcal{T}_{\alpha }=\mathcal{T}_{\alpha
}(t),t>0,$ so that we can write it as $\Psi _{n}(\mathcal{T}_{\alpha }).$
The probability law of the random time is shown to solve a fractional
diffusion equation of order $2\alpha $ and it can be expressed in terms of
Wright functions. It is interesting to stress that the introduction of a
fractional time-derivative exerts its influence only on the ``temporal''
argument, while the governing process is not affected and depends only on
the degree $n$ of the equation.

Moreover, in section 2, some particular cases of these results are analyzed:
in the non-fractional case, $\alpha =1$, we easily get $\mathcal{T}_{\alpha
}(t)\overset{a.s.}{=}t.$ For $\alpha =1/2$ it can be verified that $\mathcal{%
T}_{\alpha }$ coincides with the reflecting Brownian motion and then the
pseudoprocess governed by equation (\ref{uno.3}) reduces to $\Psi
_{n}(|B(t)|),t>0$ (where $B$ denotes a standard Brownian motion).

In section 3 the moments of $\Psi _{n}(\mathcal{T}_{\alpha })$ are obtained
in two alternative ways.

Section 4 presents some more explicit forms of the solution to equation (\ref%
{uno.3}) which can be obtained by splitting the interval of values for $%
\alpha $ in two different ones and treating them separately. Indeed, if we
restrict ourselves to the case $\alpha \in \left[ 1/2,1\right] $,\ the
distribution of the random time $\mathcal{T}_{\alpha }$ coincides with $%
\frac{1}{\alpha }\widetilde{p}_{\frac{1}{\alpha }}(u;t),$ $u\geq 0,$ where
by $\widetilde{p}_{\frac{1}{\alpha }}(\cdot ;t)$ we have denoted a stable
law of index $1/\alpha $.

As far as the other interval is concerned (i.e. $\alpha \in \left( 0,1/2%
\right] $), an explicit expression of the solution can be evaluated by
specifying $\alpha =1/m$, $m\in \mathbb{N}$, $m>2.$ In this particular
setting the pseudoprocess can be represented by $\Psi _{n}\left( G(t)\right)
,t>0,$ where $G(t)=\tprod_{j=1}^{m-1}G_{j}(t),t>0$ and $G_{j}(t),$ $%
j=1,...,m-1$ are independent random variables whose density is presented in
explicit form.

Finally we obtain some interesting results by specifying (\ref{uno.3}) for
particular values of $n.$ For example, taking $n=2$, we can conclude that
the process related to the fractional diffusion equation%
\begin{equation}
\frac{\partial ^{\alpha }}{\partial t^{\alpha }}u(x,t)=\frac{\partial ^{2}}{%
\partial x^{2}}u(x,t)\qquad x\in \mathbb{R},\;t>0,  \label{riun.12}
\end{equation}%
for $0<\alpha <1$, is represented by $B(\mathcal{T}_{\alpha })$, in
accordance with the results already known on (\ref{riun.12}). In particular
for $\alpha =1/m$, equation (\ref{riun.12}) turns out to be solved by the
density of the process $B\left( G(t)\right) ,t>0.$

In the special case $n=3$ the results above reduces to those presented in De
Gregorio (2002), while, for $n=4$, they represent a probabilistic
alternative to the analytic approach provided by Agrawal (2000).

\section{First expressions for the solution}

We start by considering the $n$-th order fractional equation and the
following corresponding initial-value problem, for $0<\alpha \leq 1,\;n\geq
2,$%
\begin{equation}
\left\{
\begin{array}{l}
\frac{\partial ^{\alpha }}{\partial t^{\alpha }}u(x,t)=k_{n}\frac{\partial
^{n}}{\partial x^{n}}u(x,t)\qquad x\in \mathbb{R},\;t>0 \\
u(x,0)=\delta (x)%
\end{array}%
\right.  \label{lapl.2}
\end{equation}%
where $k_{n}=(-1)^{q+1}$ for $n=2q,$ $q\in \mathbb{N}$, while $k_{n}=\pm 1$
for $n=2q+1$ and $\delta (\cdot )$ is the Dirac delta function. The first
step consists in evaluating the Laplace transform of the solution $u_{\alpha
}(x,t),$ namely%
\begin{equation}
U_{\alpha }(x,s)=\int_{0}^{\infty }e^{-st}u_{\alpha }(x,t)dt,  \label{lapl.1}
\end{equation}%
and recognizing that it is related to the Laplace transform of the solution $%
p_{n}(x,t)$ of the corresponding non-fractional $n$-th order equation (\ref%
{n.1}) (which can be derived from (\ref{lapl.2}) for $\alpha =1$). The
latter is usually expressed as%
\begin{equation}
p_{n}(x,t)=\frac{1}{2\pi }\int_{-\infty }^{+\infty }e^{ixz+k_{n}t(iz)^{n}}dz.
\label{sol}
\end{equation}

\

\noindent \textbf{Theorem 2.1 \ }Let $\Phi _{n}(x,s)=\int_{0}^{\infty
}e^{-st}p_{n}(x,t)dt$ be the Laplace transform of the solution to (\ref{n.1}%
); then (\ref{lapl.1}) can be expressed as follows
\begin{equation}
U_{\alpha }(x,s)=s^{\alpha -1}\Phi _{n}(x,s^{\alpha }).  \label{trasf}
\end{equation}

\noindent \textbf{Proof}

By taking the Laplace transform of (\ref{lapl.2}) and considering the
initial condition, we get%
\begin{equation}
s^{\alpha }U_{\alpha }(x,s)-s^{\alpha -1}\delta (x)=k_{n}\frac{\partial ^{n}%
}{\partial x^{n}}U_{\alpha }(x,s).  \label{lapl.4}
\end{equation}

Then, by integrating (\ref{lapl.4}) with respect to $x$ in $\left[
-\varepsilon ,\varepsilon \right] $ and letting $\varepsilon \rightarrow 0$,
we have the following condition for the $(n-1)$-th derivative%
\begin{equation*}
-s^{\alpha -1}=k_{n}\left. \frac{\partial ^{n-1}}{\partial x^{n-1}}U_{\alpha
}(x,s)\right| _{x=0^{-}}^{x=0^{+}},
\end{equation*}%
which must be added to the continuity conditions in zero holding for the $j$%
-th derivatives, for $j=0,..,n-2.$ Therefore our problem is reduced to the $%
n $-th order linear equations%
\begin{equation}
\left\{
\begin{array}{l}
k_{n}\frac{\partial ^{n}}{\partial x^{n}}U_{\alpha }(x,s)=s^{\alpha
}U_{\alpha }(x,s),\qquad x\neq 0 \\
\left. \frac{\partial ^{j}}{\partial x^{j}}U_{\alpha }(x,s)\right|
_{x=0^{-}}^{x=0^{+}}=0,\qquad \text{for }j=0,1,...,n-2 \\
\left. \frac{\partial ^{n-1}}{\partial x^{n-1}}U_{\alpha }(x,s)\right|
_{x=0^{-}}^{x=0^{+}}=-k_{n}s^{\alpha -1}%
\end{array}%
\right. .  \label{new}
\end{equation}

If we now impose the boundedness condition for $x\rightarrow \pm \infty $,
we obtain%
\begin{equation}
U_{\alpha }(x,s)=\left\{
\begin{array}{c}
\sum_{k\in I}c_{k}e^{\theta _{k}s^{\alpha /n}x},\qquad \text{if }x>0 \\
\sum_{k\in J}d_{k}e^{\theta _{k}s^{\alpha /n}x},\qquad \text{if }x\leq 0%
\end{array}%
\right. ,  \label{lapl.3}
\end{equation}%
where $\theta _{k}$ are the $n$-th roots of $k_{n},$ $I=\left\{ k:\mathbb{R}%
e(\theta _{k})<0\right\} $ and $J=\left\{ k:\mathbb{R}e(\theta
_{k})>0\right\} .$ The $n$ unknown constants $c_{k},k\in I$ and $d_{k},k\in
J,$ appearing in (\ref{lapl.3}) must be determined by taking into account
the matching conditions in (\ref{new}), as follows:%
\begin{equation}
\left\{
\begin{array}{l}
\sum_{k\in I}c_{k}\theta _{k}^{j}-\sum_{k\in J}d_{k}\theta _{k}^{j}=0,\qquad
\text{for }j=0,...,n-2 \\
\sum_{k\in I}c_{k}\theta _{k}^{n-1}-\sum_{k\in J}d_{k}\theta
_{k}^{n-1}=-k_{n}s^{\alpha /n-1}%
\end{array}%
\right. .  \label{lapl.6}
\end{equation}

By defining
\begin{equation}
z_{k}=\left\{
\begin{array}{c}
c_{k},\qquad \text{if }k\in I \\
-d_{k},\qquad \text{if }k\in J%
\end{array}%
\right. ,  \label{lapl.9}
\end{equation}%
the linear system in (\ref{lapl.6}) can be rewritten as the following
Vandermonde system%
\begin{equation}
\sum_{k=0}^{n-1}z_{k}\theta _{k}^{j}=\left\{
\begin{array}{l}
0,\qquad \text{for }j=0,...,n-2 \\
-k_{n}s^{\alpha /n-1},\qquad \text{for }j=n-1%
\end{array}%
\right. .  \label{lapl.7}
\end{equation}

Following an argument similar to Beghin and Orsingher (2005) (see p.1024-5)
we get%
\begin{eqnarray}
z_{k} &=&(-1)^{n}k_{n}s^{\alpha /n-1}\prod_{\underset{r\neq k}{r=0}}^{n-1}%
\frac{1}{\theta _{r}-\theta _{k}}  \label{lapl.8} \\
&=&\left\{
\begin{array}{l}
-\frac{1}{n}s^{\alpha /n-1}e^{\frac{2k\pi i}{n}},\qquad \text{if }k_{n}=1 \\
-\frac{1}{n}s^{\alpha /n-1}e^{\frac{(2k+1)\pi i}{n}},\qquad \text{if }%
k_{n}=-1%
\end{array}%
\right. ,  \notag
\end{eqnarray}%
where, in the last step, we have used formula (2.19) obtained therein. We
now substitute into (\ref{lapl.3}) the constants evaluated in (\ref{lapl.8}%
), taking into account (\ref{lapl.9}) and distinguishing the case of $n$
even from the odd one$.$ Indeed, for $n=2q+1,$ the roots of $k_{n}$ are
respectively%
\begin{equation}
\theta _{k}=\left\{
\begin{array}{l}
e^{\frac{2k\pi i}{n}},\qquad \text{for }k_{n}=1 \\
e^{\frac{(2k+1)\pi i}{n}},\qquad \text{for }k_{n}=-1%
\end{array}%
\right.  \label{lapl.11}
\end{equation}%
so that (\ref{lapl.3}) becomes, in this case,%
\begin{equation}
U_{\alpha }(x,s)=\left\{
\begin{array}{l}
-\frac{1}{n}s^{\alpha /n-1}\sum_{k\in I}\theta _{k}e^{\theta _{k}s^{\alpha
/n}x},\qquad \text{for }x>0 \\
\frac{1}{n}s^{\alpha /n-1}\sum_{k\in J}\theta _{k}e^{\theta _{k}s^{\alpha
/n}x},\qquad \text{for }x\leq 0%
\end{array}%
\right. .  \label{lapl.10}
\end{equation}

Analogously, for $n=2q$ and $k_{n}=(-1)^{q+1},$ the roots are $\theta
_{k}=e^{\frac{(2k+q+1)\pi i}{n}}$ so that we get%
\begin{equation}
\theta _{k}=\left\{
\begin{array}{l}
e^{\frac{(2k+q+1)\pi i}{n}}=e^{\frac{2k\pi i}{n}},\qquad \text{\ for }k_{n}=1
\\
e^{\frac{(2k+q+1)\pi i}{n}}=e^{\frac{(2k+1)\pi i}{n}},\qquad \text{\ for }%
k_{n}=-1%
\end{array}%
\right. ,  \label{lapl.12}
\end{equation}%
where, in the first line, we have used the following relationship%
\begin{equation*}
e^{(q+1)\pi i}=(-1)^{q+1}=k_{n}=1,
\end{equation*}%
while, in the second one, we have considered the fact that%
\begin{equation*}
e^{q\pi i}=(-1)k_{n}=1.
\end{equation*}

Since (\ref{lapl.12}) coincides with (\ref{lapl.11}) we obtain even for $%
n=2q $ formula (\ref{lapl.10}). The proof is completed by comparing it with
formula (12) of Lachal (2003), which reads%
\begin{equation*}
\Phi _{n}(x,s)=\left\{
\begin{array}{l}
-\frac{1}{n}s^{1/n-1}\sum_{k\in I}\theta _{k}e^{\theta _{k}s^{1/n}x},\qquad
\text{for }x>0 \\
\frac{1}{n}s^{1/n-1}\sum_{k\in J}\theta _{k}e^{\theta _{k}s^{1/n}x},\qquad
\text{for }x\leq 0%
\end{array}%
\right. .
\end{equation*}%
\hfill $\square $

\

By inverting the Laplace transform (\ref{trasf}) we can obtain a first
expression of the solution in terms of a fractional integral of a particular
stable law. Following the notation of Samorodnitsky and Taqqu (1994), we
will denote by $S_{\alpha }(\sigma ,\beta ,\mu )$ the distribution of a
\emph{stable random variable} $X$ of index $\alpha $, with characteristic
function%
\begin{equation}
Ee^{isX}=\exp \left\{ -\sigma ^{\alpha }|s|^{\alpha }\left( 1-i\beta
(sign\,s)\tan \frac{\pi \alpha }{2}\right) +i\mu s\right\} ,\qquad \alpha
\neq 1,\text{ }s\in \mathbb{R}.  \label{stable}
\end{equation}%
Moreover let $I_{(1-\alpha )}$ denote the \emph{Riemann-Liouville fractional
integral} of order $1-\alpha $, which is defined as $I_{(1-\alpha )}\left[
f(w)\right] (t)=\frac{1}{\Gamma (1-\alpha )}\int_{0}^{t}(t-w)^{-\alpha
}f(w)dw$ (see Samko et al. (1993), p.33).

\

\noindent \textbf{Theorem 2.2 \ }Let $\overline{p}_{\alpha }(\cdot ;u)$ be
the stable distribution $S_{\alpha }(\sigma ,1,0),$ with parameters $\sigma
=(u\cos \pi \alpha /2)^{1/\alpha },\beta =1,\mu =0$, then the fundamental
solution to (\ref{lapl.2}) can be expressed, for $0<\alpha <1,$ as%
\begin{equation}
u_{\alpha }(x,t)=\int_{0}^{\infty }p_{n}(x,u)I_{(1-\alpha )}\left[ \overline{%
p}_{\alpha }(w;u)\right] (t)du.  \label{riun.5}
\end{equation}

\noindent \textbf{Proof}

We recall that, for $0<\alpha \leq 2$ and $\alpha \neq 1,$ a stable random
variable $X\sim S_{\alpha }(\sigma ,1,0)$ has Laplace transform%
\begin{equation*}
E(e^{-sX})=e^{-\frac{\sigma ^{\alpha }}{\cos \left( \pi \alpha /2\right) }%
s^{\alpha }},\qquad s>0
\end{equation*}%
(see Samorodnitsky and Taqqu (1994), p.15, for details), so that, in our
case (for $\sigma =(u\cos \pi \alpha /2)^{1/\alpha }$), it reduces to $%
E(e^{-sX})=e^{-s^{\alpha }u}.$ Therefore we can rewrite (\ref{trasf}) as%
\begin{eqnarray}
U_{\alpha }(x,s) &=&s^{\alpha -1}\int_{0}^{+\infty }e^{-s^{\alpha
}t}p_{n}(x,t)dt  \label{sol.1} \\
&=&s^{\alpha -1}\int_{0}^{+\infty }\left( \int_{0}^{+\infty }e^{-sz}%
\overline{p}_{\alpha }(z;u)dz\right) p_{n}(x,u)du  \notag \\
&=&s^{\alpha -1}\int_{0}^{+\infty }e^{-sz}\left( \int_{0}^{+\infty }%
\overline{p}_{\alpha }(z;u)p_{n}(x,u)du\right) dz.  \notag
\end{eqnarray}%
For $0<\alpha <1$ the first term appearing in (\ref{sol.1}) can be easily
inverted by considering that%
\begin{equation*}
s^{\alpha -1}=\frac{1}{\Gamma (1-\alpha )}\int_{0}^{+\infty
}e^{-st}t^{-\alpha }dt
\end{equation*}%
so that \ the inverse Laplace transform of (\ref{sol.1}) can be written as%
\begin{eqnarray}
u_{\alpha }(x,t) &=&\frac{1}{\Gamma (1-\alpha )}\int_{0}^{t}(t-w)^{-\alpha
}\left( \int_{0}^{+\infty }\overline{p}_{\alpha }(w;u)p_{n}(x,u)du\right) dw
\label{sol.2} \\
&=&\frac{1}{\Gamma (1-\alpha )}\int_{0}^{+\infty }\left(
\int_{0}^{t}(t-w)^{-\alpha }\overline{p}_{\alpha }(w;u)dw\right)
p_{n}(x,u)du.  \notag
\end{eqnarray}

Finally we recognize in the last expression a fractional Riemann-Liouville
integral $I_{(1-\alpha )}$ of order $1-\alpha $ of the stable density (where
the integration is intended with respect to the first argument, since the
second represents a constant in the scale parameter). \hfill $\square $

\

The previous result suggests that the solution to our problem can be
described as the transition function $p_{n}=p_{n}(x,u)$ of a pseudoprocess $%
\Psi _{n}$ with a time-argument $\mathcal{T}_{\alpha }$ which is itself
random. Only for $\alpha =1$ we can derive from Theorem 2.1 the obvious
result that $\mathcal{T}_{\alpha }(t)\overset{a.s.}{=}t$, so that the
solution to (\ref{lapl.2}) coincides, as expected, with $p_{n}(x,t)$.\ In
all other cases the governing process coincides with the non-fractional one,
while the introduction of a fractional time-derivative exerts its influence
only on the time argument (as remarked above).

\

To check that $\mathcal{T}_{\alpha }$ possesses a true probability density
we can observe that it is non-negative: indeed it coincides with the
fractional integral of a stable density $S_{\alpha }(\sigma ,1,0)$ with
skewness parameter equal to $1$ (which, by the way, for $0<\alpha <1$, has
support restricted to $\left[ 0,\infty \right) $). Moreover it integrates to
one, as can be ascertained by the following steps:%
\begin{eqnarray*}
&&\int_{0}^{\infty }\frac{du}{\Gamma (1-\alpha )}\int_{0}^{t}(t-w)^{-\alpha }%
\overline{p}_{\alpha }(w;u)dw \\
&=&\frac{1}{\Gamma (1-\alpha )}\int_{0}^{\infty
}du\int_{0}^{t}(t-w)^{-\alpha }dw\frac{1}{2\pi i}\int_{-i\infty }^{+i\infty
}e^{sw}e^{-s^{\alpha }u}ds \\
&=&\frac{1}{2\pi i\Gamma (1-\alpha )}\int_{0}^{t}(t-w)^{-\alpha
}dw\int_{0}^{\infty }du\int_{-i\infty }^{+i\infty }e^{sw}e^{-s^{\alpha }u}ds
\\
&=&\frac{1}{2\pi i\Gamma (1-\alpha )}\int_{0}^{t}(t-w)^{-\alpha
}dw\int_{-i\infty }^{+i\infty }s^{-\alpha }e^{sw}ds \\
&=&\frac{1}{\Gamma (\alpha )\Gamma (1-\alpha )}\int_{0}^{t}w^{\alpha
-1}(t-w)^{-\alpha }dw=\frac{B(\alpha ,1-\alpha )}{\Gamma (\alpha )\Gamma
(1-\alpha )}=1\;.
\end{eqnarray*}

Our aim is now to explicit, by means of successive steps, the density $%
\overline{v}_{2\alpha }=\overline{v}_{2\alpha }(u,t)$ of $\mathcal{T}%
_{\alpha }(t),t>0$: we first prove that it satisfies a fractional diffusion
equation of order $2\alpha $ and, as a consequence, it can be expressed in
terms of \emph{Wright function}. Let
\begin{equation*}
W(x;\eta ,\beta )=\sum_{k=0}^{\infty }\frac{x^{k}}{k!\Gamma (\eta k+\beta )}
\end{equation*}%
be a Wright function of parameters $\eta ,\beta $, then we state the
following result.

\

\noindent \textbf{Theorem 2.3 \ }The fundamental solution to (\ref{lapl.2})
coincides with%
\begin{equation}
u_{\alpha }(x,t)=\int_{0}^{\infty }p_{n}(x,u)\overline{v}_{2\alpha }(u,t)du,
\label{art.3}
\end{equation}%
where
\begin{equation}
\overline{v}_{2\alpha }(u,t)=\frac{1}{t^{\alpha }}W\left( -\frac{u}{%
t^{\alpha }};-\alpha ,1-\alpha \right) ,\qquad u\geq 0,\;t>0.  \label{art.4}
\end{equation}

\noindent \textbf{Proof}

It is proved in Orsingher and Beghin (2004) that, for $0<\alpha <1$,
\begin{equation*}
I_{(1-\alpha )}\left[ \overline{p}_{\alpha }(|w|;u)\right] (t)=\frac{1}{%
\Gamma (1-\alpha )}\int_{0}^{t}(t-w)^{-\alpha }\overline{p}_{\alpha
}(|w|;u)dw
\end{equation*}%
coincides with the solution $v_{2\alpha }(u,t)$ of the following
initial-value problem, for $0<\alpha <1,$%
\begin{equation}
\left\{
\begin{array}{l}
\frac{\partial ^{2\alpha }}{\partial t^{2\alpha }}v(u,t)=\frac{\partial ^{2}%
}{\partial u^{2}}v(u,t)\qquad u\in \mathbb{R},\;t>0 \\
v(u,0)=\delta (u) \\
\frac{\partial }{\partial t}v(u,0)=0 \\
\lim_{|u|\rightarrow \infty }v(u,t)=0%
\end{array}%
\right. ,  \label{art.2}
\end{equation}%
where the second initial condition applies only for $\alpha \in \left(
1/2,1\right) .$ As a consequence, formula (\ref{riun.5}) can be rewritten as
(\ref{art.3}) with

\begin{equation}
\overline{v}_{2\alpha }(u,t)=\left\{
\begin{array}{l}
2v_{2\alpha }(u,t),\qquad \text{for }u\geq 0 \\
0,\qquad \text{for }u<0%
\end{array}%
\right. .  \label{art.1}
\end{equation}

Since it is known (see, among the others, Mainardi (1996)) that the solution
to (\ref{art.2}) can be expressed as%
\begin{eqnarray*}
v_{2\alpha }(u,t) &=&\frac{1}{2t^{\alpha }}\sum_{k=0}^{\infty }\frac{%
(-|u|t^{-\alpha })^{k}}{k!\Gamma (-\alpha k+1-\alpha )} \\
&=&\frac{1}{2t^{\alpha }}W\left( -\frac{|u|}{t^{\alpha }};-\alpha ,1-\alpha
\right) ,\qquad u\in \mathbb{R},\;t>0,
\end{eqnarray*}%
we immediately get (\ref{art.4}).\hfill$\square $

\

\noindent \textbf{Remark 2.1}

By means of the previous result we can remark again that the random time $%
\mathcal{T}_{\alpha }$ possesses a true probability density, which is
concentrated on the positive half line and moreover it is possible, thanks
to representation (\ref{art.4}), to evaluate the moments of any order $%
\delta \geq 0$ of this distribution. We recall the well known expression of
the inverse of the Gamma function as integral on the Hankel contour%
\begin{equation*}
\frac{1}{\Gamma (x)}=\frac{1}{2\pi i}\int_{Ha}e^{\tau }\tau ^{-x}d\tau ,
\end{equation*}%
which implies the representation of the Wright function as%
\begin{eqnarray*}
W(x;\eta ,\beta ) &=&\sum_{k=0}^{\infty }\frac{x^{k}}{k!\Gamma \left( \eta
k+\beta \right) } \\
&=&\frac{1}{2\pi i}\int_{Ha}e^{\tau }\sum_{k=0}^{\infty }\frac{x^{k}\tau
^{-\eta k-\beta }}{k!}d\tau \\
&=&\frac{1}{2\pi i}\int_{Ha}\tau ^{-\beta }e^{\tau +x\tau ^{-\eta }}d\tau .
\end{eqnarray*}%
Therefore we can show that%
\begin{eqnarray}
&&\int_{0}^{+\infty }u^{\delta }\overline{v}_{2\alpha }(u,t)du  \label{art.5}
\\
&=&\int_{0}^{\infty }\frac{u^{\delta }}{t^{\alpha }}W\left( -\frac{u}{%
t^{\alpha }};-\alpha ,1-\alpha \right) du  \notag \\
&=&\int_{0}^{\infty }\frac{u^{\delta }}{t^{\alpha }}\frac{du}{2\pi i}%
\int_{Ha}e^{y-\frac{u}{t^{\alpha }}y^{\alpha }}y^{\alpha -1}dy  \notag \\
&=&\frac{1}{2\pi i}\int_{Ha}e^{y}y^{\alpha -1}dy\frac{1}{t^{\alpha }}%
\int_{0}^{+\infty }e^{-\frac{u}{t^{\alpha }}y^{\alpha }}u^{\delta }du  \notag
\\
&=&\frac{t^{\alpha \delta }}{2\pi i}\int_{Ha}e^{y}y^{-\alpha \delta
-1}dy\int_{0}^{+\infty }e^{-z}z^{\delta }dz  \notag \\
&=&\frac{\Gamma (1+\delta )t^{\alpha \delta }}{\Gamma (1+\alpha \delta )}=%
\frac{t^{\alpha \delta }\Gamma (\delta )}{\alpha \Gamma (\alpha \delta )}\;.
\notag
\end{eqnarray}

From (\ref{art.5}) it is again evident that $\int_{0}^{+\infty }\overline{v}%
_{2\alpha }(u,t)du=1$ by choosing $\delta =0$.

\

\noindent \textbf{Remark 2.2}

It is interesting to analyze the particular case obtained for $\alpha =1/2$:
indeed, from the previous results, we can show that the process governed by $%
\frac{\partial ^{1/2}}{\partial t^{1/2}}u(x,t)=k_{n}\frac{\partial ^{n}}{%
\partial x^{n}}u(x,t),$ $x\in \mathbb{R},$ $t>0,$ can be represented as $%
\Psi _{n}\left( |B(t)|\right) ,t>0,$ where $B(t),t>0$ denotes a standard
Brownian motion. This can be seen by noting that $S_{1/2}\left( \frac{u^{2}}{%
2},1,0\right) $ coincides with the L\'{e}vy distribution, so that the
fractional integral in (\ref{riun.5}) reduces to%
\begin{eqnarray}
I_{(1/2)}\left[ \overline{p}_{1/2}(w;u)\right] (t) &=&\frac{1}{\Gamma (1/2)}%
\int_{0}^{t}\frac{ue^{-u^{2}/4w}}{2\sqrt{\pi (t-w)w^{3}}}dw  \label{art.6} \\
&=&\frac{e^{-u^{2}/4t}}{\sqrt{\pi t}}\;,\qquad u>0,t>0,  \notag
\end{eqnarray}%
where the second step follows by applying formula n.3.471.3, p.384 of
Gradshteyn and Rhyzik (1994), for $\mu =1/2$. Formula (\ref{art.6})
represents the density of a Brownian motion with reflecting barrier in $u=0$%
. This result is confirmed by noting that equation (\ref{art.2}), for $%
\alpha =1/2$, reduces to the heat equation $\frac{\partial }{\partial t}%
v(x,t)=\frac{\partial ^{2}}{\partial x^{2}}v(x,t)$ and then the
corresponding process coincides with a Brownian motion with $\sigma ^{2}=2t$%
. Alternatively, from (\ref{art.4}), by applying some known properties of
the Gamma function, we can write%
\begin{eqnarray}
\overline{v}_{1}(u,t) &=&\frac{1}{\sqrt{t}}\sum_{k=0}^{\infty }\frac{%
(-ut^{-1/2})^{k}}{k!\Gamma \left( 1-\frac{k+1}{2}\right) }  \label{art.7} \\
&=&\frac{1}{\sqrt{t}}\sum_{\underset{k\text{ even}}{k=0}}^{\infty }\frac{%
(-1)^{k/2}(ut^{-1/2})^{k}\Gamma \left( \frac{k+1}{2}\right) }{\pi k!}  \notag
\\
&=&\frac{1}{\pi \sqrt{t}}\sum_{\underset{k\text{ even}}{k=0}}^{\infty }\frac{%
(-1)^{k/2}(ut^{-1/2})^{k}\Gamma \left( k+1\right) \sqrt{\pi }2^{1-(k+1)}}{%
k!\Gamma \left( \frac{k}{2}+1\right) }  \notag \\
&=&\frac{1}{\sqrt{\pi t}}\sum_{j=0}^{\infty }\frac{(-1)^{j}u^{2j}(4t)^{-j}}{%
j!}=\frac{e^{-u^{2}/4t}}{\sqrt{\pi t}}\;.  \notag
\end{eqnarray}

\section{On the moments of the solution}

We are now interested in evaluating the moments of the solution of equation (%
\ref{lapl.2}), that is the moments of the pseudoprocess $\Psi _{n}(\mathcal{T%
}_{\alpha }(t)),t>0$: as we will see, they can be obtained in two
alternative ways.

By using the representation of the solution derived in (\ref{art.3}) and
thanks to the independence of the leading process from the (random) temporal
argument, we can write the $r$-th order moments as%
\begin{eqnarray}
&&E\left( \Psi _{n}^{r}(\mathcal{T}_{\alpha }(t))\right)  \label{mom.1} \\
&=&\int_{0}^{\infty }E\Psi _{n}^{r}(s)\overline{v}_{2\alpha }(s,t)ds,  \notag
\end{eqnarray}%
for $r\in \mathbb{N},t>0$. The moments of the non-fractional $n$-th order
pseudoprocess can be evaluated by means of the Fourier transform of the
solution of equation (\ref{n.1}) which can be expressed as follows%
\begin{eqnarray}
E\left( e^{i\beta \Psi _{n}(t)}\right) &=&\int_{-\infty }^{+\infty
}e^{i\beta x}p_{n}(x,t)dx=e^{(-i\beta )^{n}k_{n}t}  \label{car.1} \\
&=&\sum_{j=0}^{\infty }\frac{(i\beta )^{nj}}{(nj)!}\frac{%
(-1)^{nj}k_{n}^{j}t^{j}(nj)!}{j!}\;.  \notag
\end{eqnarray}

Therefore we get%
\begin{equation*}
E\Psi _{n}^{r}(t)=\left\{
\begin{array}{l}
\frac{(-1)^{r}(k_{n}t)^{r/n}r!}{(r/n)!}\qquad r=nj,\;j=1,2,... \\
0\qquad r\neq nj%
\end{array}%
\right. ,
\end{equation*}%
which, inserted together with (\ref{art.5}) into (\ref{mom.1}), gives, for $%
r=nj$, $j=1,2,...$%
\begin{eqnarray}
&&E\left( \Psi _{n}^{r}(\mathcal{T}_{\alpha }(t))\right)  \label{mom.5} \\
&=&\frac{(-1)^{nj}k_{n}^{j}(nj)!}{j!}\int_{0}^{\infty }s^{j}\overline{v}%
_{2\alpha }(s,t)ds  \notag \\
&=&(-1)^{nj}k_{n}^{j}t^{\alpha j}\frac{\Gamma (nj+1)}{\Gamma (\alpha j+1)}\;,
\notag
\end{eqnarray}%
while it is equal to zero for $r\neq nj$.

We can alternatively derive the moments of the pseudoprocesses by evaluating
them directly from the characteristic function of the solution. The latter
can be obtained by performing successively the Fourier and Laplace
transforms of equation (\ref{lapl.2}) as follows: let us denote by $%
\widetilde{u}_{\alpha }(\beta ,t)$ the Fourier transform of the solution,
i.e.%
\begin{equation*}
\widetilde{u}_{\alpha }(\beta ,t)=\int_{-\infty }^{+\infty }e^{i\beta
x}u_{\alpha }(x,t)dx,\qquad \beta ,t>0,
\end{equation*}%
then we get form (\ref{lapl.2})%
\begin{equation}
\frac{\partial ^{\alpha }\widetilde{u}_{\alpha }}{\partial t^{\alpha }}%
(\beta ,t)=k_{n}(-i\beta )^{n}\widetilde{u}_{\alpha }(\beta ,t).
\label{mom.2}
\end{equation}

By applying now the Laplace transform to (\ref{mom.2}) we get%
\begin{equation*}
s^{\alpha }\widetilde{U}_{\alpha }(\beta ,s)-s^{\alpha -1}=k_{n}(-i\beta
)^{n}\widetilde{U}_{\alpha }(\beta ,s),
\end{equation*}%
so that the Fourier-Laplace transform of the solution can be written as%
\begin{equation}
\widetilde{U}_{\alpha }(\beta ,s)=\frac{s^{\alpha -1}}{s^{\alpha
}-k_{n}(-i\beta )^{n}}\;.  \label{mom.3}
\end{equation}

Now recall that for the \emph{Mittag-Leffler function}
\begin{equation*}
\mathrm{E}_{\alpha ,\beta }(z)=\sum_{k=0}^{\infty }\frac{z^{k}}{\Gamma
(\alpha k+\beta )}
\end{equation*}%
the Laplace transform (for $\beta =1$) is equal to%
\begin{equation*}
\int_{0}^{\infty }e^{-sz}\mathrm{E}_{\alpha ,1}(cz^{\alpha })dz=\frac{%
s^{\alpha -1}}{s^{\alpha }-c}
\end{equation*}%
(see Podlubny (1999), formula (1.80) p. 21, for $k=0,$ $\beta =1$); hence
from (\ref{mom.3}) we get the following expression for the characteristic
function of the solution%
\begin{equation}
\widetilde{u}_{\alpha }(\beta ,t)=\mathrm{E}_{\alpha ,1}(k_{n}(-i\beta
)^{n}t^{\alpha })  \label{mom.4}
\end{equation}%
and for the solution itself
\begin{equation}
u_{\alpha }(x,t)=\frac{1}{2\pi }\int_{-\infty }^{+\infty }e^{-ix\beta }%
\mathrm{E}_{\alpha ,1}(k_{n}(-i\beta )^{n}t^{\alpha })d\beta .  \label{car.5}
\end{equation}

In the particular case $\alpha =1$ the Mittag-Leffler function reduces to
the exponential so that (\ref{mom.4}) coincides with the Fourier transform
of the solution to the $n$-th order equation, reported in (\ref{car.1}), as
it should be in the non-fractional case. Analogously, from (\ref{car.5}) we
get the usual expression of $p_{n}(x,t)$ reported in (\ref{sol}). On the
other hand, in the fractional case ($\alpha \neq 1$) formula (\ref{mom.4})
reduces, for $n=2$, to the well-known Fourier transform of the solution to
equation (\ref{riun.12}).

Finally we can evaluate the moments of the solution by rewriting formula (%
\ref{mom.4}) as%
\begin{equation*}
\widetilde{u}_{\alpha }(\beta ,t)=\sum_{j=0}^{\infty }\frac{(i\beta )^{nj}}{%
(nj)!}\frac{(-1)^{nj}k_{n}^{j}t^{\alpha j}}{\Gamma (\alpha j+1)}\Gamma
(nj+1),
\end{equation*}%
so that we get again expression (\ref{mom.5}).

\section{More explicit forms of the solution}

In order to obtain a more explicit form of the solution to (\ref{lapl.2}),
in terms of known densities, we need to distinguish between two intervals of
values for $\alpha $.

\

\noindent \textbf{(i) Case }$1/2\leq \alpha <1$

If we restrict ourselves to the case $\alpha \in \left[ 1/2,1\right) ,$ so
that $1\leq 2\alpha <2$, it is possible to apply a result obtained in Fujita
(1990), which expresses the solution to a time-fractional diffusion equation
in terms of a stable density of appropriate index. By adapting that result
to our case, we can conclude that the solution to (\ref{art.2}) coincides
with
\begin{equation*}
v_{2\alpha }(u,t)=\frac{1}{2\alpha }\widetilde{p}_{1/\alpha }(|u|;t),\qquad
u\in \mathbb{R}\text{, }t>0,
\end{equation*}%
where $\widetilde{p}_{1/\alpha }(\cdot ;t)$ denotes a stable density of
index $1/\alpha \in \left[ 1,2\right) $ with parameters $\sigma =(t\cos (\pi
-\frac{\pi }{2\alpha }))^{\alpha },$ $\beta =-1,$ $\mu =0$ (for brevity $%
S_{1/\alpha }(\sigma ,-1,0)$).

Therefore the density of $\mathcal{T}_{\alpha }(t),t>0$ is proportional to
the positive branch of a stable density, as the following expression shows:%
\begin{equation}
\overline{v}_{2\alpha }(u,t)=\frac{1}{\alpha }\widetilde{p}_{1/\alpha
}(u;t),\qquad u>0\text{, }t>0.  \label{stable2}
\end{equation}

\

\noindent \textbf{Remark 4.1}

It is possible to recognize, in the previous expression, a known density, by
resorting to results on the supremum of stable processes (see, for example,
Bingham (1973)). More precisely, let us define $Y(t)=\sup_{0\leq s\leq
t}X_{1/\alpha }(s)$ where $X_{1/\alpha }(t),t>0$ is a stable process of
index $1/\alpha $ and with characteristic function%
\begin{equation*}
E(e^{isX_{1/\alpha }(t)})=\exp \left\{ -t|s|^{1/\alpha }\left( 1+i\tan \frac{%
\pi }{2\alpha }\frac{s}{|s|}\right) \right\} ,\qquad t,s>0.
\end{equation*}

It corresponds, for any fixed $t$, to the stable law $\widetilde{p}%
_{1/\alpha }(\cdot ;t)$ defined above and, for $t$ varying, to a spectrally
negative process, which has no positive jumps (since, for $\beta =-1$, the L%
\'{e}vy-Khinchine measure assigns zero mass to $\left( 0,\infty \right) $,
see Samorodnitsky and Taqqu (1994), p.6). Under these circumstances and for $%
1/\alpha \in \left[ 1,2\right) ,$ it is known that the Laplace transform of $%
Y(t)$ is equal, for any $s,t>0,$ to%
\begin{equation*}
E(e^{-sY(t)})=\mathrm{E}_{\alpha ,1}(-st^{\alpha }),
\end{equation*}%
where $\mathrm{E}_{\alpha ,\beta }(x)$ is the Mittag-Leffler function
defined above.\emph{\ }Since it is also well-known that
\begin{equation*}
\int_{0}^{\infty }e^{-su}\widetilde{p}_{1/\alpha }(u;t)du=\alpha \mathrm{E}%
_{\alpha ,1}(-st^{\alpha }),\qquad t,s>0,
\end{equation*}%
we can conclude that
\begin{equation*}
E(e^{-sY(t)})=\int_{0}^{\infty }e^{-su}\frac{1}{\alpha }\widetilde{p}%
_{1/\alpha }(u;t)du.
\end{equation*}

Alternatively it can be shown, by adapting the result of Bingham (1973),
that the density of $Y(t)$ can be written as%
\begin{eqnarray*}
P\left\{ Y(t)\in du\right\} &=&\frac{t^{-\alpha }}{\alpha \pi }%
\sum_{n=1}^{\infty }\frac{(-1)^{n-1}}{n!}\sin \left( \pi n\alpha \right)
\Gamma \left( 1+n\alpha \right) \left( \frac{u}{t^{\alpha }}\right) ^{n-1}du
\\
&=&\frac{1}{\alpha }\widetilde{p}_{1/\alpha }(u;t)du,\qquad u>0\text{, }t>0
\end{eqnarray*}%
which coincides with (\ref{stable2}).

Formula (\ref{stable2}) shows that, for $1/2\leq \alpha <1$,%
\begin{equation*}
I_{(1-\alpha )}\left[ \overline{p}_{\alpha }(w;u)\right] (t)=\frac{1}{\alpha
}\widetilde{p}_{1/\alpha }(u;t),\qquad u>0\text{, }t>0.
\end{equation*}

Then, as a result of the fractional integration of the stable density $%
\overline{p}_{\alpha }(\cdot ;t),$ which is totally skewed to the right
(with support $\left[ 0,\infty \right) $), we obtain the positive
(normalized) branch of a new stable density $\widetilde{p}_{1/\alpha }(\cdot
;t)$ (defined on the whole real axes, since it is $1/\alpha \in \left( 1,2%
\right] $), which represents the distribution of the maximum of a stable
process of index $1/\alpha .$

\

\noindent \textbf{(ii) Case }$0<\alpha \leq 1/2$

We turn now to the other interval of values for $\alpha ,$ i.e. $\left( 0,1/2%
\right] ,$ so that, in this case, it is $0<2\alpha \leq 1$. An explicit
expression of the solution can be evaluated by specifying $\alpha =1/m$, $%
m\in \mathbb{N}$, $m>2.$ In this particular setting, problem (\ref{art.2})
becomes
\begin{equation}
\left\{
\begin{array}{l}
\frac{\partial ^{2/m}}{\partial t^{2/m}}v(u,t)=\frac{\partial ^{2}}{\partial
u^{2}}v(u,t),\qquad u\in \mathbb{R},\;t>0 \\
v(u,0)=\delta (u) \\
\lim_{|u|\rightarrow \infty }v(u,t)=0%
\end{array}%
\right.  \label{fcaa}
\end{equation}%
so that it can be considered as a special case of the fractional telegraph
equation studied in Beghin and Orsingher (2003), for $\lambda =0$ and $c=1$.
By applying formula (2.11) of that paper, the solution to (\ref{fcaa}) can
be expressed, for $u\in \mathbb{R},$ $t>0,$ as%
\begin{eqnarray}
&&v_{2/m}(u,t)  \label{vu2m} \\
&=&\left( \frac{m}{2\pi }\right) ^{\frac{m-1}{2}}\frac{1}{2\sqrt{t}}%
\int_{0}^{\infty }dw_{1}...\int_{0}^{\infty }dw_{m-1}\cdot  \notag \\
&&\cdot e^{-\frac{w_{1}^{m}+...+w_{m-1}^{m}}{\sqrt[m-1]{m^{m}t}}}w_{2}\cdot
\cdot \cdot w_{m-1}^{m-2}\left[ \delta (u-w_{1}\cdot \cdot \cdot
w_{m-1})+\delta (u+w_{1}\cdot \cdot \cdot w_{m-1}\right] .  \notag
\end{eqnarray}

By taking, as before,%
\begin{equation}
\overline{v}_{2/m}(u,t)=\left\{
\begin{array}{l}
2v_{2/m}(u,t),\qquad \text{for }u\geq 0 \\
0,\qquad \text{for }u<0%
\end{array}%
\right. ,  \label{fcaa1}
\end{equation}
the solution to our problem (\ref{lapl.2}) can be expressed, in this case, as%
\begin{eqnarray*}
u_{1/m}(x,t) &=&\int_{0}^{\infty }p_{n}(x,u)\overline{v}_{2/m}(u,t)du \\
&=&\int_{0}^{\infty }p_{n}(x,u)p_{G(t)}(u)du
\end{eqnarray*}%
where $G(t)=\tprod_{j=1}^{m-1}G_{j}(t),t>0$ and $G_{j}(t),$ $j=1,...,m-1$
are independent random variables with the following probability law
\begin{equation}
p_{G_{j}(t)}(w)=\frac{1}{m^{\frac{j}{m-1}-1}t^{\frac{j}{m(m-1)}}\Gamma (%
\frac{j}{m})}\exp \left( -\frac{w^{m}}{\sqrt[m-1]{m^{m}t}}\right)
w^{j-1}\quad w>0.  \label{dens}
\end{equation}

We can check the independence by noting that%
\begin{eqnarray}
&&\tprod_{j=1}^{m-1}p_{G_{j}(t)}(w_{j})  \label{dens.2} \\
&=&\tprod_{j=1}^{m-1}\frac{1}{m^{\frac{j}{m-1}-1}t^{\frac{j}{m(m-1)}}\Gamma (%
\frac{j}{m})}\exp \left( -\frac{w_{j}^{m}}{\sqrt[m-1]{m^{m}t}}\right)
w_{j}^{j-1}  \notag \\
&=&\left( \frac{m}{2\pi }\right) ^{\frac{m-1}{2}}\frac{1}{\sqrt{t}}\exp
\left( -\frac{\sum_{j=1}^{m-1}w_{j}^{m}}{\sqrt[m-1]{m^{m}t}}\right)
\tprod_{j=1}^{m-1}w_{j}^{j-1},  \notag
\end{eqnarray}%
where, in the second step, we have applied the multiplication formula of the
Gamma function%
\begin{equation*}
\Gamma (z)\Gamma \left( z+\frac{1}{m}\right) ...\Gamma \left( z+\frac{m-1}{m}%
\right) =\left( 2\pi \right) ^{\frac{m-1}{2}}m^{\frac{1}{2}-mz}\Gamma \left(
mz\right) ,
\end{equation*}%
for $z=1/m$. The last expression in (\ref{dens.2}) coincides with the joint
density of the variables $G_{j}(t)$ given in formula (1.7) of Beghin and
Orsingher (2003). Therefore the corresponding pseudoprocess is represented,
in this case, as $\Psi _{n}\left( G(t)\right) ,t>0.$

\ \newpage

\noindent \textbf{Remark 4.2}

We can check the previous results, obtained separately for the two
intervals, by choosing $\alpha =1/2$. From both cases we obtain again that
the pseudoprocess governed by our equation can be represented by $\Psi
_{n}\left( |B(t)|\right) ,t>0.$

Indeed from the first case, i.e. for\emph{\ }$1/2\leq \alpha <1,$ we get, by
means of (\ref{stable2}), that the density of $\mathcal{T}_{\alpha }(t),t>0,$
for $\alpha =1/2$, coincides with the folded normal. More precisely, $S_{2}(%
\sqrt{t},-1,0)$ coincides with $N(0,2t)$ and then%
\begin{equation}
\overline{v}_{1}(u,t)=2\widetilde{p}_{2}(u;t)=\frac{e^{-u^{2}/4t}}{\sqrt{\pi
t}}  \label{v1}
\end{equation}%
for $u>0$, $t>0.$

On the other hand, if we consider the expression of the density of $\mathcal{%
T}_{\alpha }$ obtained for $0<\alpha \leq 1/2$, we get, for $\alpha =1/2$
and $m=2$, from (\ref{dens.2}) that again
\begin{equation}
p_{G_{1}(t)}(u)=\frac{e^{-u^{2}/4t}}{\sqrt{\pi t}}\text{.}  \label{v2}
\end{equation}

Moreover both (\ref{v1}) and (\ref{v2}) coincide with (\ref{art.6}) derived
above, as expected.

\

An interesting application of our results can be obtained by specializing
Theorems 2.2 and 2.3 to the particular case $n=2.$ In this situation the
pseudoprocess $\Psi _{n}(t),t>0$ reduces to the Brownian motion (with
variance $2t$) $B(t),t>0$ and therefore the solution of (\ref{lapl.2})
coincides with the transition density of the process $B(T_{\alpha }(t)),t>0$
obtained by the composition of $B$ with the random time $T_{\alpha }.$ We
state this last result as follows

\

\noindent \textbf{Corollary 4.1}

The solution of the problem%
\begin{equation}
\left\{
\begin{array}{l}
\frac{\partial ^{\alpha }}{\partial t^{\alpha }}u(x,t)=\frac{\partial ^{2}}{%
\partial x^{2}}u(x,t)\qquad x\in \mathbb{R},\;t>0, \\
u(x,0)=\delta (x)%
\end{array}%
\right.  \label{rem}
\end{equation}%
for $0<\alpha \leq 1,$ is represented by the transition function of $B(%
\mathcal{T}_{\alpha }).$ The density of the random time $\mathcal{T}_{\alpha
}=\mathcal{T}_{\alpha }(t),t>0$ is the folded solution of the
time-fractional equation%
\begin{equation*}
\frac{\partial ^{2\alpha }}{\partial t^{2\alpha }}v(u,t)=\frac{\partial ^{2}%
}{\partial u^{2}}v(u,t)\qquad u\in \mathbb{R},\;t>0
\end{equation*}%
and is given in (\ref{art.4}).

\

We can prove that this is in accordance with what is already known on (\ref%
{rem}): for $n=2$ we can substitute in (\ref{art.3}) the transition function
of the Brownian motion, so that we get:

\begin{eqnarray}
u_{\alpha }(x,t) &=&\frac{1}{t^{\alpha }}\int_{0}^{\infty }\frac{%
e^{-x^{2}/4u}du}{\sqrt{4\pi u}}W\left( -\frac{u}{t^{\alpha }};-\alpha
,1-\alpha \right)  \label{rem.2} \\
&=&\frac{1}{t^{\alpha }}\int_{0}^{\infty }\frac{e^{-x^{2}/4u}du}{\sqrt{4\pi u%
}}\frac{1}{2\pi i}\int_{Ha}\frac{e^{y-\frac{u}{t^{\alpha }}y^{\alpha }}}{%
y^{1-\alpha }}dy  \notag \\
&=&\frac{1}{4it^{\alpha }\sqrt{\pi ^{3}}}\int_{Ha}\frac{e^{y}}{y^{1-\alpha }}%
dy\int_{0}^{\infty }\frac{e^{-\frac{x^{2}}{4u}-\frac{u}{t^{\alpha }}%
y^{\alpha }}}{\sqrt{u}}du.  \notag
\end{eqnarray}

If we prove now that
\begin{equation}
\int_{0}^{\infty }\frac{e^{-\frac{x^{2}}{4u}-\frac{u}{t^{\alpha }}y^{\alpha
}}}{\sqrt{u}}du=\sqrt{\pi }t^{\alpha /2}y^{-\alpha /2}e^{-\frac{|x|}{%
t^{\alpha /2}}y^{\alpha /2}}  \label{rem.3}
\end{equation}%
and substitute (\ref{rem.3}) into (\ref{rem.2}), we finally get the known
result%
\begin{eqnarray*}
u_{\alpha }(x,t) &=&\frac{1}{2t^{\alpha /2}}\frac{1}{2\pi i}\int_{Ha}\frac{%
e^{y-\frac{|x|}{t^{\alpha /2}}y^{\alpha /2}}}{y^{1-\alpha /2}}dy \\
&=&\frac{1}{2t^{\alpha /2}}W\left( -\frac{|x|}{t^{\alpha /2}};-\frac{\alpha
}{2},1-\frac{\alpha }{2}\right) .
\end{eqnarray*}

In order to verify formula (\ref{rem.3}) we use the following relationship,
known for the Laplace transform of the first-passage time of the Brownian
motion,%
\begin{equation*}
e^{-|x|\sqrt{s}}=\int_{0}^{\infty }e^{-su}\frac{|x|}{2\sqrt{\pi }\sqrt{u^{3}}%
}e^{-\frac{|x|^{2}}{4u}}du,
\end{equation*}%
which, integrated with respect to $x$ gives (\ref{rem.3}), for $s=y^{\alpha
}/t^{\alpha }.$ Alternatively, we can apply formula n.3.471.9, p.384 of
Gradshteyn and Ryzhik (1994), for $\beta =x^{2}/4$, $\gamma =y^{\alpha
}/t^{\alpha },$ $\nu =1/2$ (noting that $K_{1/2}(z)=\sqrt{\pi /2z}e^{-z}$,
see Gradshteyn and Ryzhik (1994), n.8469.3, p.978).

\

If we restrict ourselves to the case $\alpha \in \left[ 1/2,1\right) ,$ the
density of $\mathcal{T}_{\alpha }(t),t>0$ is again proportional to the
positive branch of a stable density, as expressed in (\ref{stable2}).

On the other hand, for $\alpha \in \left( 0,1/2\right] $ and in particular $%
\alpha =1/m,\;m\in \mathbb{N}$, it is represented by the law presented in (%
\ref{vu2m}) and (\ref{fcaa1}), so that the process governed by equation (\ref%
{rem}) is, in this case, $B(\prod_{j=1}^{m-1}G_{j}(t)),t>0.$

\

\begin{center}
{\Large References}
\end{center}

\noindent \textbf{Anh V.V.,Leonenko N.N. (2000)},\textbf{\ ``}Scaling laws
for fractional diffusion-wave equation with singular data'', \emph{%
Statistics and Probability Letters, }\textbf{48},\textbf{\ }239-252.

\noindent \textbf{Agrawal O.P. (2000)}, ``A general solution for the
fourth-order fractional diffusion-wave equation'', \emph{Fract. Calc. Appl.
Anal.}, \textbf{3}, 1, 1-12.

\noindent \textbf{Angulo J.M., Ruiz-Medina M.D., Anh V.V., Grecksch W. (2000)%
},\textbf{\ }Fractional diffusion and fractional heat equation, \emph{%
Advances in Applied Probability,} \textbf{32} (4), 1077-1099.

\noindent \textbf{Beghin L., Orsingher E. (2003), }``The telegraph process
stopped at stable-distributed times and its connection with the fractional
telegraph equation'', \emph{Fract. Calc. Appl. Anal.}, \textbf{6} (2),
187-204.

\noindent \textbf{Beghin L., Orsingher E. (2005), }``The distribution of the
local time for `pseudoprocesses' and its connection with fractional
diffusion equations'', \emph{Stoch. Proc. Appl.}, \textbf{115}, 1017--1040.

\noindent \textbf{Beghin L., Orsingher E., Ragozina T. (2001), }``Joint
Distributions of the Maximum and the Process for Higher-Order Diffusions'',
\emph{Stoch. Proc. Appl.}, \textbf{94}, 71-93.

\noindent \textbf{Beghin L., Hochberg K., Orsingher E. (2000), }%
``Conditional Maximal Distributions of Processes Related to Higher-Order
Heat-Type Equations'', \emph{Stoch. Proc. Appl.}, \textbf{85}, 209-223.

\noindent \textbf{Bingham N.H. (1973)},\textbf{\ }``Maxima of sums of random
variables and suprema of stable processes'', \emph{%
Z.Wahrscheinlichkeitstheorie verw. Geb., }\textbf{26},\textbf{\ }273-296.

\noindent \textbf{Daletsky Yu. L. (1969),} ``Integration in function
spaces'', In: \emph{Progress in Mathematics}$,$ R.V. Gamkrelidze, ed.,
\textbf{4}, 87-132.

\noindent \textbf{De Gregorio A. (2002)}, \textquotedblleft Pseudoprocessi
governati da equazioni frazionarie di ordine superiore al
secondo\textquotedblright , \emph{Tesi di Laurea}, Universit\`{a} di Roma
\textquotedblleft La Sapienza\textquotedblright .

\noindent \textbf{Fujita Y. (1990), }Integrodifferential equation which
interpolates the heat equation and the wave equation (I), \emph{Osaka
Journal of Mathematics, }\textbf{27}, 309-321.

\noindent \textbf{Gradshteyn I.S., Rhyzik I.M. (1994)}\emph{, Tables of
Integrals, Series and Products, }Alan Jeffrey Editor, Academic Press, London.

\noindent \textbf{Hochberg K. J., Orsingher E. (1994),} ``The arc-sine law
and its analogs for processes governed by signed and complex measures'',
\emph{Stoch. Proc. Appl., }\textbf{52}, 273-292.

\noindent \textbf{Hochberg K. J., Orsingher E. (1996),} ``Composition of
stochastic processes governed by higher-order parabolic and hyperbolic
equations'', \emph{Journal of Theoretical Probability, }\textbf{9} (2),
511-532.

\noindent \textbf{Krylov V. Yu. (1960),} ``Some properties of the
distribution corresponding to the equation $\frac{\partial u}{\partial t}%
=(-1)^{p+1}\frac{\partial ^{2p}u}{\partial x^{2p}}$'', \emph{Soviet Math.
Dokl., }\textbf{1}, 260-263.

\noindent \textbf{Lachal A. (2003), }``Distributions of sojourn time,
maximum and minimum for pseudo-processes governed by higher-order heat-type
equations'', \emph{Electronic Journ. Prob., }\textbf{8} (20), 1-53.

\noindent \textbf{Lachal A. (2007), }``First hitting time and place for
pseudo-processes driven by the equation $\frac{\partial }{\partial t}=\pm
\frac{\partial ^{N}}{\partial x^{N}}$ subject to a linear drift'', \emph{%
Stoch. Proc. Appl.}.

\noindent \textbf{Mainardi F. (1996)}, ``Fractional relaxation-oscillation
and fractional diffusion-wave phenomena'', \emph{Chaos, Solitons and Fractals%
}, \textbf{7}, 1461-1477.

\noindent \textbf{Nigmatullin R.R. (1986)}, ``The realization of the
generalized transfer equation in a medium with fractal geometry'', \emph{%
Phys. Stat. Sol.}, (b) \textbf{133}, 425-430.

\noindent \textbf{Nikitin Y., Orsingher E. (2000),} ``Conditional sojourn
distribution of `processes' related to some higher-order heat-type
equations'', \emph{Journ. Theor. Probability}, \textbf{13} (4), 997-1012.

\noindent \textbf{Nishioka K. (1997),} ``The first hitting time and place of
a half-line by a biharmonic pseudo process'' \emph{Japanese Journal of
Mathematics, }\textbf{23} (2), 235-280.

\noindent \textbf{Orsingher E. (1991),} ``Processes governed by signed
measures connected with third-order 'heat-type' equations'', \emph{Lith.
Math. Journ.}, \textbf{31}, 321-334.

\noindent \textbf{Orsingher E., Beghin L. (2004),} ``Time-fractional
equations and telegraph processes with Brownian time'', \emph{Probability
Theory and Related Fields.}, \textbf{128}, 141-160.

\noindent \textbf{Podlubny I. (1999)}, \emph{Fractional Differential
Equations}, Academic Press, San Diego, (1999).

\noindent \textbf{Saichev A.I., Zaslavsky G.M. (1997)}, ``Fractional kinetic
equations: solutions and applications'', \emph{Chaos}, \textbf{7} (4),
753-764.

\noindent \textbf{Samko S.G., Kilbas A.A., Marichev O.I. (1993)}, \emph{%
Fractional Integrals and Derivatives: Theory and Applications}, Gordon and
Breach Science Publishers.

\noindent \textbf{Samorodnitsky G., Taqqu M.S. (1994)}, \emph{Stable
Non-Gaussian Random Processes}, Chapman and Hall, New York.

\noindent \textbf{Schneider W.R., Wyss W. (1989)}, ``Fractional diffusion
and wave equations'', \emph{Journal of Mathematical Physics}, \textbf{30},
134-144.\newline
\textbf{Wyss W. (1986)}, ``Fractional diffusion equation'', \emph{Journal of
Mathematical Physics}, \textbf{27}, 2782-2785.\newline
\textbf{Wyss W. (2000)}, ``The fractional Black-Scholes equation'', \emph{%
Fractional Calculus and Applied Analysis}, \textbf{3} (1), 51-61.\newline

\end{document}